\documentclass[journal]{IEEEtran}

\usepackage[utf8]{inputenc}
\usepackage{graphicx}
\usepackage{cite}       
\graphicspath{ {Figures_bnd_wave/} }
\usepackage{amsmath}
\usepackage{amsfonts}
\usepackage{color}
\usepackage{amssymb}
\usepackage{mathrsfs}
\usepackage{float}
\usepackage{enumerate}
\usepackage{verbatim}
\usepackage{setspace}
\usepackage{epsfig}
\usepackage{url}
\usepackage{graphicx}
\usepackage{lscape}

\newtheorem{innercustomassumption}{Assumption}
\newenvironment{customassumption}[1]
{\renewcommand\theinnercustomassumption{#1}\innercustomassumption}
{\endinnercustomassumption}

\newtheorem{innercustomrem}{Remark}
\newenvironment{customrem}[1]
{\renewcommand\theinnercustomrem{#1}\innercustomrem}
{\endinnercustomrem}

    \newtheorem{innercustomprob}{Problem}
\newenvironment{customprob}[1]
{\renewcommand\theinnercustomprob{#1}\innercustomprob}
{\endinnercustomprob}

\newtheorem{innercustomthm}{Theorem}
\newenvironment{customthm}[1]
{\renewcommand\theinnercustomthm{#1}\innercustomthm}
{\endinnercustomthm}

\begin{document}

\title{Internal stabilization of three interconnected semilinear reaction-diffusion PDEs with one actuated state}		
%

\author{Constantinos~Kitsos, Rami~Katz and
Emilia~Fridman,
\thanks{C. Kitsos ({\tt\small constantinos.kitsos@ec-nantes.fr}) is with the Institut de Recherche en Génie Civil et Mécanique (UMR CNRS 6183), Ecole Centrale de Nantes. R. Katz ({\tt\small ramikatz@mail.tau.ac.il}) and E. Fridman ({\tt\small  emilia@eng.tau.ac.il}) are with the School of Electrical Engineering, Tel Aviv University, Israel. }
\thanks{Supported by  Israel Science Foundation (grant no. 673/19) and
by Chana and Heinrich Manderman Chair at Tel Aviv University.}%
}
\markboth{}%
{}

\maketitle

\begin{abstract}
This work deals with the  exponential stabilization of a  system of three semilinear parabolic partial differential equations (PDEs), written in a strict feedforward form. The diffusion coefficients are considered distinct and the PDEs are interconnected via both a reaction matrix and a nonlinearity. Only one of the PDEs is assumed to be controlled internally, thereby leading to an underactuated system. Constructive and efficient control of such underactuated systems is a nontrivial open problem, which has been solved recently for the linear case. In this work, these results are extended to the semilinear case, which is  highly challenging due the interconnection that is introduced by the nonlinearity. Modal decomposition is employed, where due to nonlinearity, the finite-dimensional part of the solution is coupled with the infinite-dimensional tail. A transformation is then performed to map the finite-dimensional part into a target system, which allows for an efficient design of a static linear proportional state-feedback controller. Furthermore, a high-gain approach is employed in order to compensate for the nonlilinear terms. Lyapunov stability analysis is performed, leading to LMI conditions guaranteeing exponential stability with arbitrary decay rate. The LMIs are shown to always be feasible, provided the number of actuators and the value of the high gain parameter are large enough. Numerical examples show the efficiency of the proposed approach.   
\end{abstract}
\section{Introduction}
Interconnected reaction-diffusion equations are widely used in science and engineering, including models in population dynamics and spatial ecology, chemical reactions, magnetic systems and evolutionary game theory (see e.g. \cite{kopell1973plane,okubo2001diffusion}, and the surveys \cite{ruiz2022control}, \cite[Section 3.3]{pierre2010global}). Notably, 
tumor growth models \cite{gallinato2017tumor}, diffusional neural networks \cite{chen2017intermittent} and nuclear reactors \cite{jiang2014semilinear} are characterized by parabolic PDEs with nonlinear interconnections as the ones we consider here and limited number of controlled states (\emph{underactuated systems}). 

For scalar semilinear reaction-diffusion equations, constructive methods for output-feedback stabilization via distributed or boundary finite-dimensional controller are given in \cite{katz2020constructive}. Prior to that, stabilization via LQR-based controllers was achieved for a scalar heat equation in \cite{Hagen03}. Controllability of a scalar semilinear heat equation was studied in \cite{fabre1995approximate} and further results on stabilization of nonlinear PDEs appeared in \cite{christofides1998robust,coron2004global}.

In comparison to the scalar case, control of interconnected and underactuated reaction-diffusion systems of equations introduces additional theoretical challenges. For example, the presence of distinct diffusion coefficients often leads to significant technical difficulties in analyzing the system and deriving a suitable control law (see e.g. survey \cite{ammar2011recent} and the controllability results therein). Moreover, from an implementation perspective, high computational costs are often required to compute the controller gains and implement the controller. The problem of control of underactuated systems was  first introduced in \cite{lions1989remarques} (see also \cite{kitsosFirdman22b} for a more exhaustive literature review). Stabilization of underactuated systems of linear reaction-diffusion PDEs with multiple states and a single controlled equation was presented in \cite{KitsosFridman22a,kitsosFirdman22b}, for the first time, for both distributed and boundary control. 

The present work extends the results of \cite{KitsosFridman22a,kitsosFirdman22b} to PDE systems interconnected by \emph{unknown} Lipschitz nonlinearities. We consider a system written in a strict feedforward form, whose diffusion matrix is diagonal and contains distinct elements. The remaining terms consist of a reaction matrix and a triangular Lipschitz nonlinearity. We assume that only the first PDE is controlled, with $N$ actuators acting internally on the state. To design a stabilizing control law, we follow a modal decomposition approach, (see e.g. \cite{katz2020constructive,katz2023global} for scalar semilinear PDEs, and \cite{kitsosFirdman22b} for PDE systems) to obtain a finite-dimensional part and an infinite-dimensional tail, which are coupled due to the nonlinearity. Then, the challenge in stabilizing the system becomes twofold: on the one hand, the underactuation and distinct diffusion coefficients complicate the efficient and scalable controller design, and on the other hand, the nonlinearity that couples all the modes, thereby causing a spillover behavior \cite{Hagen03}. To deal with the first challenge, we employ a transformation, recently introduced in \cite{KitsosFridman22a}, to transform the finite-dimensional part into an equivalent target system. The obtained target system allows for an efficient and scalable state-feedback controller design, with significantly reduced computational complexity. The infinite-dimensional analog of this transformation first appeared in \cite{KitsosPhd,kitsos2020high} to deal with a more general reaction and under stronger regularity assumptions on the system solutions. Here, modal decomposition allows less demanding requirements of the system. To overcome the second challenge, we exploit the triangularity of the system and combine the controller design with a high-gain approach inspired by finite-dimensional systems (see \cite{khalil2017high} for instance about high-gain observers) and initiated first for infinite-dimensional systems in  \cite{kitsos2020high}. We perform Lyapunov stability analysis of the closed-loop system, leading to LMIs. Differently from finite-dimensional systems \cite{khalil2017high}, here the high-gain parameter also intervenes in the infinite-dimension tail, thereby leading to a much more complicated stability analysis. Given an arbitrary decay rate, the LMIs are shown to always be feasible provided the number of actuators $N$ and the high-gain parameter are large enough. The derived LMIs provide a constructive and scalable stabilization algorithm.

The paper is organized as follows. The problem formulation and assumptions are presented in Section \ref{sec:problem}. Our internal stabilization approach is presented in Section \ref{sec:internal}, where Theorem \ref{theorem1} constitutes the main result of this paper. Section \ref{sec:example} presents a numerical example. Conclusions are drawn in Section \ref{sec:conclusion}. 

\textit{Notations:} For $x\in {{\mathbb{R}}^{m}}$, $\left| x \right|$ denotes its usual Euclidean norm. For a matrix $Q\in {{\mathbb{R}}^{m\times m}}$, ${Q}^{\top}$ denotes its transpose, $\left| Q \right|:=\sup \left\{ \left| Qw \right|,\left| w \right|=1 \right\}$ is its induced norm, $\mathrm{Sym}(Q)=\frac{Q+Q^{\top}}{2}$ stands for its symmetric part and $\sigma_{\min}(Q)$, $\sigma_{\max}(Q)$ denote its minimum and maximum eigenvalues, respectively. By $\text{diag}\{ A_1, \ldots, A_m \}$ we denote the (block) diagonal matrix with diagonal elements $A_1,\ldots,A_m$. $I_m$ denotes the identity matrix of dimension $m$ and $\otimes$ denotes the Kronecker product. For $f, g$ in $L^2\left( 0,L;\mathbb R^m\right)$, we denote their inner product by $\left <f,g \right >=\int_0^L f^\top (x)g(x) dx$ with induced norm $\Vert \cdot \Vert_{L^2}$ ($L^2\left( 0,L;\mathbb R^m\right)$ is the space of equivalence classes of measurable functions $f:[0,L]\to\mathbb{R}^m$ with $\Vert f \Vert_{L^2}<+\infty$).  By $\ell^2(\mathbb N;\mathbb R^m)$ we denote the Hilbert space of square summable sequences $x=(x_n)_{n=1}^{+\infty}\subseteq \mathbb{R}^m$. 

\section{Problem formulation}\label{sec:problem}

Consider the following system of three coupled parabolic PDEs, where for $(t,x)\in [0,+\infty)\times(0,L)$:
\begin{subequations}\label{sys}
	\begin{align}
		&z_t(t,x)=D z_{xx}(t,x)+Q z(t,x)+ f\left (z(t,x)\right )\notag\\&\hspace{14mm}+B \sum_{j=1}^{N}b_j(x)u_j(t),\label{sys0}\\
		\label{BC} &\gamma_{11} z(t,0)+\gamma_{12} z_x(t,0)=0, \notag\\&\gamma_{21} z(t,L)+\gamma_{22} z_x(t,L)=0,\\
		\label{IC} &z(0,x)=z^0(x).
	\end{align}
\end{subequations}
Here, the state is given by $z=\operatorname{col} \{z_1, z_2,  z_3\}$, the diffusion matrix $$D=
\begin{pmatrix}
	d_{1} &0 &0 \\
	0& d_2 &0 \\
	0& 0& d_{3}
\end{pmatrix}$$ consists of possibly distinct coefficients $d_1,d_2,d_3>0$ and   $Q$ is a reaction matrix split as follows: 
\begin{align} Q&= Q_0+ Q_1;\notag\\
	Q_0&={
		\begin{pmatrix}
			0& 0&0\cr 
			q_{2,1}
			&0&0\\0&q_{3,2}&0
		\end{pmatrix}, \quad \ Q_1=\begin{pmatrix}
			q_{1,1}& q_{1,2}&q_{1,3}\cr 
			0
			&q_{2,2}&q_{2,3}\\0&0&q_{3,3}
		\end{pmatrix}\label{eq:Qmatdef}
	}.
\end{align}
The nonlinear term $f$ admits a triangular form as following:
\begin{gather} \label{f}
	f(z)=\begin{pmatrix}f_1(z_1,z_2,z_3)\\ f_2(z_2,z_3)\\ f_3(z_3) \end{pmatrix}
\end{gather}
and $B = \begin{pmatrix}
	1\\0\\0
\end{pmatrix}$ is the control coefficient. The boundary conditions are defined using constants $\left\{\gamma_{ij}\right\}_{i,j=1}^2\subseteq \mathbb R$ satisfying $\gamma_{i1}^2+\gamma_{i2}^2\neq 0,\quad  i=1,2.$  

The system \eqref{sys} represents an \emph{underactuated} system in which the scalar control inputs $\left\{u_j\right\}_{j=1}^N$, with $N$ to be determined below, are multiplied by the shape functions $\left\{b_j\right\}_{j=1}^N\subseteq L^2(0,L)$ and act internally on the first equation only. The control input then propagates to the equations for the states $z_2$ and $z_3$ via the reaction matrix Q. The nonlinear term (in particular, $f_1$) then interconnects $z_2$ and $z_3$ back to $z_1$.

		
		\begin{customassumption}{1}\label{assumption}\it 
			The elements of $Q_0$ are nonzero. Namely, $q_{2,1}, q_{3,2} \neq 0$.
		\end{customassumption}
		Assumption \ref{assumption} guarantees controllability of $(Q_0,B)$.
		
		\begin{customassumption}{2}\label{assumptionF}\it 
			The function $f$ is in ${Lip}\left (L^2\left (0,L;\mathbb R^3\right ) \right )$. More precisely, there exists constants $\ell_i > 0, \quad i=1,2,3$ such that the following inequalities hold:
			\begin{align}
				\Vert f_i(z)-f_i(\bar z) \Vert_{L^2} \leq \ell_i\Vert z-\bar z \Vert_{L^2} \label{Lipschitz}
			\end{align}
			for every $z, \bar z\in L^2\left (0,L;\mathbb R^3\right ).$ Moreover, $f_i(0)\equiv 0, \ i=1,2,3.$
		\end{customassumption}
		The above global Lipschitz assumption is selected here in order to simplify the already cumbersome derivations. However, it can be weakened to a locally Lipschitz property, which is more common in reaction-diffusion systems. In this case, the approach presented in this manuscript allows to obtain a local stabilization result, with explicit estimate on the domain of attraction (see, e.g., \cite{Katz2022} for similar arguments).
		\begin{customrem}{1}\label{rem:f2=f3=0}\it
			We will henceforth assume that $f_2\equiv f_3 \equiv 0$, whereas $f_1$ is nonzero with no condition on its Lipschitz constant $\ell_1$. Subject to this assumption we will design a controller and derive LMI conditions which guarantee $L^2$ stability of the closed-loop system. Then, by robustness of LMIs to small perturbations, the proposed control strategy will also work with nonzero $f_2$ and $f_3$, provided that $\ell_2$ and $\ell_3$ are small enough. In future work, we will further provide explicit estimates on $\ell_2$ and $\ell_3$, which preserve the $L^2$-stability of the closed-loop system. Note that the interconnection of $z_1$, $z_2$ and $z_3$ holds under the assumption that only $f_1$ is nonzero.
		\end{customrem}

Consider the following family of scalar Stürm-Liouville eigenvalue problems for $i=1,2,3$:
\begin{align}
	\begin{aligned}    
		d_i  \varphi^{\prime \prime}(x)+\bar \lambda \varphi(x) =&0, \quad 0<x<L,\\
		\gamma_{11} \varphi(0)+\gamma_{12}\varphi^\prime(0)=&\gamma_{21} \varphi(L)+\gamma_{22}\varphi^\prime(L)=0, \label{SLprob}
	\end{aligned}
\end{align}
admitting a monotone sequence of simple eigenvalues $\bar \lambda_{n,i}=d_i  \lambda_n$, where $\lambda_n\propto n^2$ (see for instance \cite[Section 3.3.1]{orlov2020nonsmooth}) are the eigenvalues of \eqref{SLprob} with $d_i=1$. The eigenvalues correspond to a sequence of eigenfunctions $(\varphi_n)_{n=1}^{+\infty}$. The eigenfunctions form a complete and orthonormal system in $L^2(0,L)$. For the case of Neumann or Dirichlet boundary conditions in \eqref{sys} (i.e., when one of the pairs $(\gamma_{1,i},\gamma_{2,j}),i,j=1,2$ is zero) the eigenfunctions can be computed explicitly. For Robin boundary conditions, our approach can still be employed, by numerically computing the eigenvalues and eigenfunctions and taking into account estimates on them (see,  e.g., \cite[Section 3.3.1]{orlov2020nonsmooth}).

\begin{customassumption}{3}\label{assumptiononB} \it
	Recall the shape functions $\left\{b_j\right\}_{j=1}^N$. For any $N\in \mathbb{N}$, the matrix 
	\begin{align}
		\mathcal B_{ N\times N}:=\operatorname{col} \left\{\mathcal B_{j}^{\top} \right\}_{j=1}^{ N},\label{eq:BNNMat}
	\end{align}
	with 
	$
	\mathcal B_n:=\begin{pmatrix}
		b_{1,n}&\cdots&b_{N,n}
	\end{pmatrix}^\top$
	and $b_{j,n}:=\left\{\left<b_j,\varphi_n \right>\right\}_{j,n=1}^N$ 
	is non-singular.  Furthermore, for all $N \in \mathbb N,$ we have
	\begin{align} \label{Bboundedness}
		\left(\sum_{k=1}^N\left[\left\|b_k \right\|^2_{L^2} - \left| \mathcal{B}_k\right|^2 \right]  \right) \lvert \mathcal B_{N\times N}^{-1} \rvert^2\leq \eta \lambda_{N+1}^{\beta} \end{align}
	for some constants $\eta>0$ and $\beta \in [0,1).$
\end{customassumption}
\begin{customrem}{2}\it
	The inverse of  $\mathcal B_{ N\times N}$ will be used explicitly in the controller design. Similar assumption appears in several works on stabilization of scalar parabolic PDEs (see for instance \cite{Hagen03}). The condition \eqref{Bboundedness}, however, is a technical one, and will be used below to prove feasibility guarantees for the derived LMIs. General conditions on 
	$\left\{b_j\right\}_{j=1}^N$, which guarantee that \eqref{Bboundedness} holds, 
	do not appear in the literature to the best of our knowledge and their derivation is a highly nontrivial problem, which we leave for future work.
	Note that the conditions of Assumption \ref{assumptiononB} hold for the particular case $b_n(\cdot)= \varphi_n(\cdot)$ for all $n\in \mathbb{N}$.
\end{customrem}

We are now ready to describe the problem that we address in this paper:
\begin{customprob}{1} \label{Problem} \it
	Let Assumptions \ref{assumption}-\ref{assumptiononB} hold. Given any decay rate $\delta>0$, find the number of actuators $N$ and static linear proportional state-feedback control laws $\left\{u_i(t)\right\}_{i=1}^N$ such that for any initial condition $z^0 \in H^1\left (0,L;\mathbb R^3\right)$ satisfying the boundary conditions \eqref{BC}, the solution to the closed-loop system obtained from \eqref{sys}  satisfies
	\begin{align}
		\Vert z(t,\cdot)\Vert_{L^2} \leq M e^{-\delta t} \Vert z^0(\cdot) \Vert_{L^2}, \quad t \geq 0 \label{eq:stabilityineq}
	\end{align}
	with some $M> 0$. 
\end{customprob}
\section{Stabilization method} \label{sec:internal}

Employing modal decomposition and \eqref{SLprob}, we begin by presenting the solution of \eqref{sys} as 
\begin{align}
	z_i(t,\cdot)=\sum_{n=1}^\infty z_{i,n}(t) \varphi_n(\cdot), \quad i=1,2,3 \label{eq:Zseriespres}
\end{align}
with the projection coefficients $z_{i,n}$ given by 
\begin{gather}
	z_{i,n}=\left <z_i,\varphi_n\right >\label{zi}.
\end{gather}
Taking the time-derivative of \eqref{zi}, substituting the dynamics \eqref{sys}, and integrating by parts, we obtain the following ODEs for $z_n=\operatorname{col} \{z_{1,n}, z_{2,n},  z_{3,n} \}$:  
\begin{align}
	&\dot z_n(t)=\int_0^L z_t(t,x) \varphi_n(x)  dx\notag\\&=\left [D z_x(\cdot)\varphi_n(\cdot)-Dz(\cdot)  \varphi_n^{\prime} (\cdot)\right ]_0^L \notag\\&+\left (-\lambda_n D +Q\right )z_n(t)+F_n[z(t)]\notag\\&+B \sum_{j=1}^N u_j(t)  \int_0^L \varphi_n(x)  b_j(x) dx, \notag
\end{align}
where $F_n[z]:=\int_0^L f\left (z(t,x)\right)\phi_n(x)dx$, which by virtue of the homogeneous boundary conditions for $\varphi_n(x)$ and $z(t,x)$, can be written as follows:
\begin{align}
	\dot z_n(t)=&\left (-\lambda_n D +Q\right )z_n(t) +F_n[z(t)]+B \sum_{j=1}^{N} b_{j,n} u_j(t). \label{zn_dot0}
\end{align}
The representation \eqref{eq:Zseriespres}, as well as the derivation of \eqref{zn_dot0} are justified by existence and uniqueness of classical solution of the closed-loop system derived after \eqref{eq:ClosedLoopFin2} below.

\begin{customrem}{3}\it
	Differently from the linear case presented in \cite{KitsosFridman22a}, \eqref{zn_dot0} shows that we cannot decouple the stabilization of the slower finite-dimensional part of the solution (which corresponds to smaller eigenvalues) from the faster infinite-dimensional tail, due to the introduction of the nonlinear terms $F_n(t)$, which couple all solution modes (this phenomenon is also known as spillover, see \cite{Hagen03}).
\end{customrem}
Using the notation $z^N=\text{col}\{ z_j\}_{j=1}^N\in \mathbb R^{3N}$, we obtain the following ODE system for the finite-dimensional part:
\begin{align}
	\dot z^N(t)=A z^N(t)+F^N[z(t)]+\tilde B u(t),  \label{z^Nopen}
\end{align}
where $F^N[z]:=\text{col}\{ F_j[z]\}_{j=1}^N, u(t):=\text{col}\{ u_j(t)\}_{j=1}^N,$
\begin{align} 
	A:=&\text{diag}\{-\lambda_1 D+ Q,\ldots,-\lambda_N D+ Q \}, \label{Amatrix}
\end{align}
and  $\tilde B\in \mathbb R^{3N\times N}$ is given by $$
\tilde B:= \operatorname{col}\{ B \mathcal B_1,\ldots, B\mathcal B_{N}\}=\left ( I_N \otimes B \right )\mathcal B_{N \times N}. $$
By invoking the Hautus lemma, it is easy to see that the pair $(A, \tilde B)$ is stabilizable under Assumption \ref{assumptiononB}. In light of \eqref{z^Nopen} and \eqref{Amatrix}, a straightforward control strategy would be to seek controller gains $K\in \mathbb{R}^{1\times 3N}$ such that $\tilde{A}+\tilde{B}K$ is Hurwitz. Unfortunately, this approach is highly computationally expensive for large $N$. Instead, we propose to exploit the triangular structure of the system and the distinct diffusion coefficients by applying a suitable transformation to the finite-dimensional part \eqref{z^Nopen}. This type of transformation was first introduced in \cite[Chapter 3]{KitsosPhd,kitsos2020high} in a more general setting and in \cite{KitsosFridman22a} to solve the stabilization problem of $m$ coupled linear PDEs with one controlled state via modal decomposition. The transformation therein was subject to generalized Sylvester equation, whose solution was given explicitly. The transformation greatly simplifies the design of controller gains by reducing it to LMIs of dimension $3$, instead of $3N$ (the size of the finite-dimensional part). This is achieved by stabilizing the component of the reaction matrix $Q_0\in \mathbb{R}^{3\times 3}$, in contrast to the straightforward approach presented above. More details about such transformations appear in \cite{kitsosFirdman22b}. 

We, therefore, apply the transformation
\begin{equation}\label{eq:ModeTransf}
	y_n=T_n z_n,\quad n\geq 1    
\end{equation}
to system \eqref{zn_dot0} 
for $T_n \in \mathbb R^{3\times 3}$ given by
\begin{align}
	T_n=\left\{ \begin{array}{ll} 
		I_3+\lambda_n \begin{pmatrix}
			0 & \kappa& 0\\0&0&0\\0&0&0
		\end{pmatrix}, \quad & 1\leq  n \leq N,\\I_3, \quad &n\geq N+1\end{array}\right., \label{T_n}
\end{align}
where 
\begin{equation}\label{eq:KappaDef}
	\kappa:=\frac{d_3-d_2}{q_{2,1}}.   
\end{equation}
Note that $\mathcal T:=(T_n)_{n=1}^{+\infty}:\ell^2(\mathbb N;\mathbb R^3)\to \ell^2(\mathbb N;\mathbb R^3)$ is a bounded operator with bounded inverse. The inverse of $T_n$ is given by 
\begin{align}
	T_n^{-1}=\left\{ \begin{array}{ll} 
		I_3-\lambda_n \begin{pmatrix}
			0 & \kappa& 0\\0&0&0\\0&0&0
		\end{pmatrix}, \quad & 1\leq  n \leq N,\\I_3, \quad &n\geq N+1\end{array}\right..
\end{align}
Moreover,
{\begin{align} \label{Tn-1bound}
		\max_{n\in\mathbb N}\vert T_n^{-1}\vert,\max_{n\in\mathbb N}\vert T_n \vert \leq \sigma_N:=1+\vert \kappa \vert \lambda_N.
\end{align} }
Notice that such a transformation first appeared in \cite{KitsosFridman22a} for stabilization of a \emph{cascade} system of $m$ linear parabolic PDEs. In that work, it was also shown that a transformation $T_n$ of the triangular form \eqref{T_n} was more intricate when the number of distinct diffusion coefficients is large. For the case $m=3$
and $d_2=d_3$, $T_n$ is reduced to the identity matrix. Therefore, here we are mostly interested by the case when $d_2\neq d_3$.

Following the above transformation and recalling \eqref{eq:Qmatdef}, we write the new system (for both the finite-dimensional part and the tail) as follows:
\begin{align}  \left\{ \begin{array}{ll}  \begin{aligned}
			\dot y_n(t)&=\left (-\lambda_n d_3 I_3 + Q_0 +B G_n + J_n \right)y_n(t) + B \mathcal{B}_n^{\top}\times\\
			&{ \operatorname{col}\left\{u_j \right\}_{j=1}^N}+T_n F_n\left [\sum_{k=1}^{+\infty}T_k^{-1} y_k(t)\varphi_k\right ],  \ n \leq N,\\
			\dot y_n(t)&=\left (-\lambda_n D+  Q\right )y_n(t) + F_n\left [\sum_{k=1}^{+\infty}T_k^{-1} y_k(t)\varphi_k\right ]\\&+{ B \mathcal{B}_n^{\top}\operatorname{col}\left\{u_j \right\}_{j=1}^N},  \  n\geq N+1
	\end{aligned}\end{array}\right.  \label{target_sys0} 
\end{align}
with $\mathcal{B}_n$ given in \eqref{eq:BNNMat} and
\begin{align}
	J_n &= T_nQ_1T_n^{-1},\nonumber \\
	G_n &= \begin{pmatrix}
		\lambda_n(d_2-d_1)  & \quad \lambda_n^2\kappa\left ( d_1-d_3\right )   & \quad 0
	\end{pmatrix}.  \label{Gn}
\end{align}
To obtain a closed-loop system description for $t\geq 0$, we introduce the notations
\begin{align}
	&y^{N}(t) = \operatorname{col}\left\{ y_j(t)\right\}_{j=1}^N\in \mathbb R^{3N}, \ \mathcal{J}_N = \operatorname{diag}\left\{J_j \right\}_{j=1}^N,\nonumber \\
	&\Lambda = \operatorname{diag}\left\{\lambda_j \right\}_{j=1}^N, \ \mathcal{G}_N = \operatorname{diag}\left\{G_j \right\}_{j=1}^N,\label{eq:Fj} \\
	&F_j(t) = F_j\left [\sum_{k=1}^{+\infty}T_k^{-1} y_k(t)\varphi_k\right ], \ j\geq 1 \nonumber 
\end{align}
and present \eqref{target_sys0} as follows:
\begin{align}
	& \dot y^N(t)=\left (- d_3\Lambda \otimes I_3 + I_N \otimes Q_0 +\left( I_N \otimes B \right) \mathcal{G}_N \right)y^N(t) \nonumber \\
	&\hspace{12mm}+\mathcal J_Ny^N(t)  +\left( I_N \otimes B \right)\mathcal{B}_{N\times N} \operatorname{col}\left\{u_j \right\}_{j=1}^N \nonumber\\
	&\hspace{12mm}+\operatorname{col}\left\{T_j F_j(t)\right\}_{j=1}^N,\label{eq:PartsSeparation}\\
	&\dot y_n(t)=\left (-\lambda_n D+  Q\right )y_n(t) + F_n(t) \nonumber \\
	&\hspace{10mm}+B \mathcal{B}_n^{\top}\operatorname{col}\left\{u_j \right\}_{j=1}^N,  \quad  n\geq N+1. \nonumber
\end{align}
Since $(Q_0,B)$ is stabilizable, we further choose $K_0 \in \mathbb R^{1\times 3}$ such that $Q_0+B K_0$ is Hurwitz. We then propose the control law of the form
\begin{align}
	&\operatorname{col}\left\{u_j \right\}_{j=1}^N = \mathcal{B}_{N\times N}^{-1} \left[-\mathcal{G}_N+K \right]y^{N}(t), \nonumber\\
	&K = \operatorname{diag}\left\{K_j \right\}_{j=1}^N \in \mathbb{R}^{N\times 3N}, \label{eq:ContLaw}
\end{align}
with $K_j\in \mathbb{R}^{1\times 3}, \ 1\leq j\leq N$ to be specified shortly. Substituting \eqref{eq:ContLaw} into \eqref{eq:PartsSeparation}, we obtain for the closed-loop system
\begin{align}
	& \dot y^N(t)=\left (- d_3\Lambda \otimes I_3 + I_N \otimes  Q_0 +\left( I_N \otimes B \right) K \right)y^N(t) \nonumber \\
	&\hspace{10mm} +\mathcal{J}_Ny^N(t)+\operatorname{col}\left\{T_j F_j(t)\right\}_{j=1}^N, \label{eq:PartsSeparation1}\\
	&\dot y_n(t)=\left (-\lambda_n D+  Q\right )y_n(t) + F_n(t)\nonumber\\
	&\hspace{10mm}+B \mathcal{B}_n^{\top}\mathcal{B}_{N\times N}^{-1} \left[-\mathcal{G}_N+K \right]y^{N}(t),  \quad  n\geq N+1.\nonumber
\end{align}
Then, noting that 
\begin{align*}
	& I_N \otimes Q_0 +\left( I_N \otimes B \right) K = \operatorname{diag}\left\{ Q_0+BK_j\right\}_{j=1}^N
\end{align*}
and taking into account that $Q_0+BK_0$ is Hurwitz, we finally choose gains $K_j$ as follows:
\begin{align}
	&K_j  = \gamma^4 K_0 \Gamma^{-1};\quad \Gamma = \operatorname{diag}\left\{\gamma^{3},\gamma^{2},\gamma \right\}\nonumber \\
	&\Rightarrow K = \gamma^4 \left(I_N \otimes \left[K_0\Gamma^{-1} \right] \right).\label{eq:GammaDef}
\end{align}
Here $\gamma\geq 1$ is a high-gain tuning parameter to be determined later. Substituting \eqref{eq:GammaDef} into \eqref{eq:PartsSeparation1}, we rewrite the closed-loop system for all $t\geq 0$ as follows:
\begin{align*}
	& \dot y^N(t)=\left (- d_3\Lambda \otimes I_3 + I_N\otimes \left[ Q_0+\gamma^4B K_0 \Gamma^{-1} \right] \right)y^N(t) \\
	&\hspace{10mm}+ \mathcal{J}_Ny^N(t)+\operatorname{col}\left\{T_j F_j\right\}_{j=1}^N,\\
	&\dot y_n(t)=\left (-\lambda_n D+  Q\right )y_n(t) + F_n+B \mathcal{B}_n^{\top}\mathcal{B}_{N\times N}^{-1}\times \\
	&\hspace{10mm} \left(-\mathcal{G}_N+ \gamma^4 \left(I_N \otimes \left[K_0\Gamma^{-1} \right] \right) \right)y^{N}(t),  \quad  n\geq N+1.
\end{align*}
\begin{customrem}{4} \it
	The introduction of a high-gain controller in \eqref{eq:GammaDef} is motivated by high-gain control of nonlinear finite-dimensional systems, where stabilization of triangular systems is based on combining 
	stabilization of the linear part of the system with high-gains, which compensate the nonlinear term. Differently from finite dimensional systems, in our case the high-gain parameter also appears in the infinite-dimensional tail, both explicitly and through the nonlinearity  (see \eqref{eq:CSSemilinear}). This leads to a more challenging stability analysis and a more delicate treatment of the choice of $\gamma$ in proving feasibility guarantees for the derived LMIs. Note also, that the high-gain approach was firstly extended to infinite dimension in \cite{kitsos2020high} in the context of observers for coupled PDEs.
\end{customrem}

Finally, we introduce the following change of variables for the finite-dimensional part:
\begin{align}
	&\bar{F}_j(t) = \Gamma^{-1}T_jF_j(t), \quad  \bar{y}_j=\Gamma^{-1} y_j, \quad  \overline{J}_n =\Gamma^{-1}J_n\Gamma ,\nonumber \\
	&\bar{y}^N(t) = \operatorname{col}\left\{\bar{y}_j \right\}_{j=1}^N = (I_N\otimes \Gamma^{-1})y^{N}(t),\nonumber \\
	&\bar{F}^N(t) = \operatorname{col}\left\{\bar{F}_j(t) \right\}_{j=1}^N,\quad  \overline{\mathcal{J}}_{N}=\operatorname{diag}\left\{\overline{J}_j \right\}_{j=1}^N. \label{eq:BarTransf} 
\end{align}
Recalling Remark 1, we find that
\begin{align}
	&\overline{F}_j(t) = \gamma^{-3}F_j(t), \quad 1\leq j \leq N,\nonumber\\
	&\Rightarrow \overline{F}^N(t) = \gamma^{-3}F^{N}(t).\label{eq:OnlyF1}
\end{align}
Hence, we have the following ODEs for $\bar{y}_n(t), \ 1\leq n \leq N$ and $t\geq 0$:
\begin{align*}
	&\dot{\bar{y}}_n(t) = \Gamma^{-1}\left[-d_3 \lambda_n I_3+  { Q_0}+\gamma^4B K_0 \Gamma^{-1}\right]\Gamma \bar{y}_n(t)\\
	&\hspace{12mm}+{ \overline{J}_ny_n(t)}+\bar{F}_n(t)\\
	&\hspace{-1mm} = \left[-d_3 \lambda_n I_3 + \gamma  Q_0+\gamma^4 \Gamma^{-1}B K_0+{\overline{J}_n} \right]\bar{y}_n(t)+\gamma^{-3}F_n(t)\\
	&\hspace{-1mm} = \left[-d_3 \lambda_n I_3 + \gamma\left( Q_0+B K_0\right)+{\overline{J}_n} \right]\bar{y}_n(t)+\gamma^{-3}F_n(t).
\end{align*}
Taking the latter into account, we end up with the following closed-loop system for $t\geq 0$:
\begin{align}
	& \dot{\bar{y}}^N(t)=\left (- d_3\Lambda \otimes I_3 + \gamma I_N\otimes \left[{ Q_0}+B K_0 \right] \right)\bar{y}^N(t) \nonumber \\
	&\hspace{10mm}+{\overline{\mathcal{J}}_{N}\bar{y}^N(t)}+\gamma^{-3}F^N(t), \label{eq:ClosedLoopFin}\\
	&\dot y_n(t)=\left (-\lambda_n D+  Q\right )y_n(t) + F_n+B \mathcal{B}_n^{\top}\mathcal{B}_{N\times N}^{-1}\times \nonumber \\
	&\hspace{10mm} \left(-\mathcal{G}_N+ \gamma^4 \left(I_N \otimes \left[K_0\Gamma^{-1} \right] \right) \right)\nonumber \\
	&\hspace{10mm}\times (I_N\otimes \Gamma) \bar{y}^{N}(t),  \quad  n\geq N+1.\nonumber 
\end{align}
{Recalling \eqref{eq:Qmatdef} and $J_n$ in \eqref{Gn}, it can be verified that for $1\leq n \leq N$
	\begin{align}
		&\overline{J}_n =  \begin{pmatrix}
			q_{1,1} & \frac{q_{1,2}+\lambda_n\kappa(q_{2,2}-q_{1,1})}{\gamma} & \frac{q_{1,3}+\lambda_n \kappa q_{2,3}}{\gamma^2}\\
			0 & q_{2,2} & \frac{q_{2,3}}{\gamma}\\
			0 & 0 & q_{3,3}
		\end{pmatrix}  \nonumber \\
		&\Rightarrow \left|\overline{\mathcal{J}}_{N}\right| = \max_{1\leq n\leq N}\left|\overline{J}_n \right| \leq \max_{1\leq i\leq 3}\left|q_{i,i} \right|+ O\left(\frac{\lambda_N}{\gamma}\right).\label{eq:Jestimate}
\end{align}}
Set
\begin{align}
	\overline{\mathcal{G}}_N = \mathcal{G}_N(I_N\otimes \Gamma) = \operatorname{diag}\left\{G_j\Gamma \right\}_{j=1}^N \in \mathbb{R}^{N\times 3N}.\label{eq:CalGNDef}
\end{align}
It can be easily verified that \eqref{Gn} implies
\begin{align}
	&\gamma^{-4}\left|\overline{\mathcal{G}}_N\right|=\gamma^{-4} \max_{1\leq j \leq N}\left|G_j\Gamma \right|\leq \xi_{N,\gamma}, \label{eq:ClaGNBound}\\
	&\xi_{N,\gamma} = \max \left(\frac{\lambda_N\left| d_2-d_1\right|}{\gamma},\left(\frac{\lambda_N}{\gamma} \right)^2 \vert \kappa( d_1-d_3)\vert \right )\nonumber
\end{align}
Using \eqref{eq:ClosedLoopFin} and \eqref{eq:CalGNDef}, we write the closed-loop system as 
\begin{align}
	& \dot{\bar{y}}^N(t)=\left (- d_3\Lambda \otimes I_3 + \gamma I_N\otimes \left[{Q_0}+B K_0 \right] \right)\bar{y}^N(t) \nonumber \\
	&\hspace{10mm}{ +\overline{\mathcal{J}}_{N}\bar{y}^N(t)}+\gamma^{-3}F^N(t), \label{eq:ClosedLoopFin2}\\
	&\dot y_n(t)=\left (-\lambda_n D+  Q\right )y_n(t) + F_n(t)+B \mathcal{B}_n^{\top}\mathcal{B}_{N\times N}^{-1}\times \nonumber \\
	&\hspace{10mm} \left(-\overline{\mathcal{G}}_N+ \gamma^4 I_N \otimes K_0 \right) \bar{y}^{N}(t),  \quad  n\geq N+1.\nonumber 
\end{align}

Prior to the $L^2$ stability analysis, taking into account the Lipschitz property of $f$ and assuming $z^0(\cdot)\in H^1\left (0,L;\mathbb R^3\right )$, we may invoke classical arguments, see for instance \cite{pazy1983semigroups,pierre2010global} on the existence of unique solutions $z(\cdot,\cdot)$ of the closed-loop system \eqref{sys} with control law \eqref{eq:ContLaw}, \eqref{eq:GammaDef}, belonging to $C^1\left ((0,T];L^2\left (0,L;\mathbb R^3\right ) \right )$ for some $T>0$, with $z(\cdot,x)$ continuous on $t=0$. The stability analysis below, combined with a standard bootstrap argument then guarantees that for the solution of the closed-loop system, subject to a stabilizing controller, we have $T=+\infty$.

For $L^2$-stability analysis of the closed-loop system \eqref{eq:ClosedLoopFin}, we introduce the Lyapunov functional $\mathcal V:L^2\left (0,L;\mathbb R^3\right )\to \mathbb R_{\geq 0}$ acting on $y(t,\cdot)=\sum_{n=1}^{+\infty} y_n(t)\phi_n(\cdot)$ (with $y_j=\Gamma \bar y_j, j=1,\ldots,N$) and defined as:
\begin{equation}
	\mathcal V[y]:=\frac{1}{2}  ( \bar{y}^N)^\top \bar P  \bar{y}^N+\frac{\rho}{2} \sum_{n=N+1}^{+\infty} \vert y_n\vert^2,\label{eq:Vfunc}
\end{equation} 
where $0<\rho$ and $\bar{P} = I_N\otimes P$ with $0\prec  P\in \mathbb{R}^{3\times 3}$. Note that thanks to the transformation $T_n$ in \eqref{T_n}, each of the blocks $-d_3 \lambda _n I_3+\gamma (Q_0+BK_0), n=1,\ldots,N$ in closed-loop system \eqref{eq:ClosedLoopFin} will be stabilized simultaneously by use of the Lyapunov matrix $P$. Indeed, $P$ stabilizes the matrix $Q_0$ via the gain $K_0$, whereas $-d_3 \lambda_nI_3 P \prec 0$ for all $n$. In particular, $P$ is \emph{independent} of $N$.

We define $V(t):=\mathcal V[y](t), \forall t \in [0,T).$
Then, there exist constants $\underline{c}, \overline{c}>0$ such that 
\begin{align}
	\underline{c} \left\|z(t,\cdot) \right\|^2_{L^2} \leq {V}(t) \leq \overline{c} \left\|z(t,\cdot) \right\|^2_{L^2}, \quad  t\in [0,T). \label{eq:Comp}
\end{align}
Indeed, recalling $\Gamma$ in \eqref{eq:GammaDef}, where $\gamma\geq 1$, we have 
\begin{align*}
	&{V}(t)\leq \frac{\sigma_{max}(P)}{2}\left|\bar{y}^N(t) \right|^2+ \frac{\rho}{2}\sum_{n=N+1}^{\infty} \left|y_n(t) \right|^2 \\
	&\overset{\eqref{eq:BarTransf}}{\leq } \max \left( \frac{\sigma_{max}(P)}{2\gamma^2},\frac{\rho}{2}\right)\sum_{n=1}^{\infty} \left|y_n(t) \right|^2\\
	&\leq \max \left( \frac{\sigma_{max}(P)}{2\gamma^2},\frac{\rho}{2}\right)\sigma_N^2 \sum_{n=1}^{\infty} \left|z_n(t) \right|^2
\end{align*}
where the last inequality follows from \eqref{eq:ModeTransf}-\eqref{Tn-1bound}. Setting $\overline{c} = \max \left( \frac{\sigma_{max}(P)}{2\gamma^2},\frac{\rho}{2}\right)\sigma_N^2$ and applying Parseval's equality, we obtain the upper bound in \eqref{eq:Comp}. By similar arguments, one can determine the lower bound $\underline{c} = \min \left( \frac{\sigma_{min}(P)}{2\gamma^6},\frac{\rho}{2}\right)\sigma_N^{-2}$.

\begin{customrem}{5} \it
	Note that the derived constants $\underline{c}, \overline{c}$ and the Lyapunov analysis below imply that $M$ in \eqref{eq:stabilityineq} depends on powers of the high-gain parameter $\gamma$. Therefore, a large $\gamma$ may lead to significant overshoot in the transient behavior of the system. This phenomenon is well-known in high-gain controller design.
\end{customrem}

Differentiating ${V}(t)$ in $(0,T)$ along the solution to the closed-loop system \eqref{eq:ClosedLoopFin2}, we get
\begin{align}
	&\dot{{V}}(t) = (\bar y^N(t))^\top \left [ -d_3 \Lambda \otimes P+{\bar{P}\overline{\mathcal{J}}_{N}}+\gamma I_N\otimes\right.\notag\\
	&\left. \operatorname{Sym}\left (P\left [{ Q_0}+B K_0\right ]\right )\right ]\bar y^N(t)+\gamma^{-3}(\bar y^N(t))^\top \bar P F^N(t)\notag\\
	&+\rho\sum_{n=N+1}^{+\infty}y_n^\top(t) \left( -\lambda_n D +  \text{Sym}\left ( Q\right )\right) y_n(t)\notag\\&+\rho \sum_{n=N+1}^{+\infty}y_n^\top(t) \left[F_n(t)+\gamma^4B \mathcal B_n^\top \mathcal Q_{N\times N} \bar{y}^{N}(t) \right],
	\label{Lyapder0}
\end{align}
where 
\begin{align}
	\mathcal Q_{N\times N,\gamma}:=\mathcal B_{N\times N}^{-1}\left(-\gamma^{-4}\overline{\mathcal{G}}_N + I_N \otimes K_0 \right).
\end{align}
satisfies 
\begin{align} \label{Qnnbound}
	&\left|\mathcal Q_{N\times N,\gamma} \right|\overset{\eqref{eq:ClaGNBound}}{\leq}\left|\mathcal B_{N\times N}^{-1} \right|\left[\xi_{N,\gamma}+ \left| K_0\right| \right].
\end{align}
Let $\alpha_0>0$. By Young's inequality, we have
\begin{align*}
	&\rho \sum_{n=N+1}^{\infty}y_n^{\top}(t)F_n(t)\leq \frac{\rho}{2\alpha_0}\sum_{n=N+1}^{\infty}\left| y_n(t) \right|^2 \\
	&- \frac{\alpha_0 \rho}{2} \left| F^N(t) \right|^2+\frac{\alpha_0 \rho}{2}  \sum_{n=1}^{\infty}\left| F_n(t) \right|^2
\end{align*}
We further apply Parseval's equality to obtain 
\begin{align} 
	& \sum_{n=1}^{\infty}\left|F_n(t) \right|^2  = \int_0^L \left| F\left[\sum_{n=1}^{+\infty}T_n^{-1} y_n(t)\phi_n(x)\right]\right|^2dx\notag\\
	&\overset{\eqref{Lipschitz},\eqref{eq:Fj}}{\leq} \ell_1\int_{0}^L \left|\sum_{n=1}^{+\infty}T_n^{-1} y_n(t)\phi_n(x) \right|^2dx \notag\\
	&=\ell_1 \sum_{n=1}^{\infty}\left|T_n^{-1} y_n(t)\right|^2 \overset{\eqref{eq:BarTransf}}{=}\ell_1 \sum_{n=1}^{N}\left|T_n^{-1}\Gamma \overline{y}_n(t)\right|^2\nonumber \\
	&+\ell_1 \sum_{n=N+1}^{\infty}\left|T_n^{-1} y_n(t)\right|^2\overset{\eqref{eq:BarTransf}}{=} \ell_1 \left( \bar{y}^N(t)\right)^{\top}\mathcal{M}_{N,\gamma} \bar{y}^N(t)\notag\\
	&+\ell_1 \sum_{n=N+1}^{\infty} \left|y_n(t) \right|^2,\label{eq:CSSemilinear}
\end{align}
where $\mathcal{M}_{N,\gamma}$ is the symmetric matrix
\begin{align}
	&\mathcal{M}_{N,\gamma} = \operatorname{diag}\left\{\Gamma T_j^{-\top}T_{j}^{-1} \Gamma \right\}_{j=1}^N \prec \gamma^6\sigma_N^2 I_{3N}, \label{eq:CalMNDef}
\end{align}
where the above bound easily follows from \eqref{Tn-1bound}.
Similarly, for some $\alpha_1>0$,
\begin{align}
	&\gamma^4 \sum_{n=N+1}^{\infty}y_n^{\top}(t)B \mathcal B_n^\top \mathcal Q_{N\times N,\gamma} \bar{y}^{N}(t)\leq \frac{1}{2\alpha_1}\sum_{n=N+1}^{\infty}\left| y_n(t)\right|^2\nonumber \\
	&+\frac{\alpha_1 \gamma^8}{2}\sum_{n=N+1}^{\infty}\left|B \mathcal B_n^\top \mathcal Q_{N\times N,\gamma} \bar{y}^{N}(t) \right|^2,
\end{align}
where
\begin{align}
	&\sum_{n=N+1}^{\infty}\left|B \mathcal B_n^\top \mathcal Q_{N\times N,\gamma} \bar{y}^{N}(t) \right|^2 = \left(Q_{N\times N,\gamma} \bar{y}^{N}(t) \right)^{\top}\times \nonumber \\
	&\left(\sum_{n=N+1}^{\infty}\mathcal{B}_nB^{\top}B\mathcal{B}_n^{\top} \right)\left(Q_{N\times N,\gamma} \bar{y}^{N}(t)\right)\nonumber \\
	&\leq \left(\sum_{m=1}^N\sum_{N+1}^{\infty} \left| b_{m,n}\right|^2\right)\left|Q_{N\times N,\gamma} \bar{y}^{N}(t) \right|^2\nonumber \\
	&=\left(\sum_{k=1}^N\left[\left\|b_k \right\|_{L^2}^2 - \left| \mathcal{B}_k\right|^2 \right]  \right)\left|Q_{N\times N,\gamma} \bar{y}^{N}(t) \right|^2.\label{eq:BBnBound}
\end{align}
Finally, let $\delta>0$ be a desired decay rate and 
\begin{equation*}
	\eta(t) = \operatorname{col}\left\{\bar{y}^N(t), F^N(t) \right\}.
\end{equation*}
From \eqref{Lyapder0}-\eqref{eq:BBnBound} we find that
\begin{align}
	\dot{{V}}(t)+2\delta {V}(t) &\leq     \eta(t)^{\top} \Phi \eta(t)\\
	&+\rho \sum_{n=N+1}^{\infty}y_n^{\top}(t) W_ny_n(t)\leq 0\label{eq:Vdot}
\end{align}
provided
\begin{align}
	&W_n : =  -\lambda_n D +  \text{Sym}\left ( Q\right )+\delta I_3 \nonumber \\
	&\hspace{10mm}+ \left[\frac{1}{2\alpha_0}+\frac{\alpha_0 \ell_1}{2}+\frac{1}{2\alpha_1}\right]I_3 \prec 0, \ n>N \nonumber \\
	& \Phi = \begin{bmatrix}
		\phi & \gamma^{-3}\bar{P}\\
		* & -\frac{\alpha_0 \rho}{2}I_{3N}
	\end{bmatrix}\prec 0,\label{eq:LMI1}
\end{align}
where $\bar{P} = I_N\otimes P$ and
\begin{align*}
	& \phi = -d_3 \Lambda \otimes P+\gamma I_N\otimes \operatorname{Sym}\left (P\left [{Q_0}+B K_0+\frac{\delta}{\gamma}I_3\right ]\right )\\
	&\hspace{2mm} +\frac{\alpha_0 \rho \ell_1}{2}\mathcal{M}_{N,\gamma}+{\bar{P} \overline{\mathcal{J}}_{N}} \\
	&\hspace{2mm} +\frac{\alpha_1\rho \gamma^8}{2}\left(\sum_{k=1}^N\left[\left\|b_k \right\|_{L^2}^2 - \left| \mathcal{B}_k\right|^2 \right]  \right) \mathcal{Q}_{N\times N,\gamma}^{\top}\mathcal{Q}_{N\times N,\gamma}.
\end{align*}

By monotonicity of $\lambda_n$, we have $W_n\leq W_{N+1},$ for all $\ n\geq N+1$. By Schur complement, $W_{N+1}\prec 0$ iff
\begin{equation}\label{eq:TailSchur}
	\begin{bmatrix}
		-\lambda_{N+1} D +  \text{Sym}\left ( Q\right ) + \frac{\alpha_0 \ell_1}{2}I_3 & I_3 & I_3\\
		* & -2\alpha_0I_3 & 0\\
		* & * & -2\alpha_1I_3
	\end{bmatrix}\prec 0.
\end{equation}
Summarizing, we arrive at the main result of this work:
\begin{customthm}{1}\label{theorem1}
	Consider system \eqref{sys} with $f_2 \equiv f_3\equiv0$, subject to assumptions \ref{assumption}-\ref{assumptiononB} and the control law \eqref{eq:ContLaw}, \eqref{eq:GammaDef}, where $K_0$ is such that {$Q_0+BK_0$} is Hurwitz, $\gamma>1$ is a high-gain tuning parameter and  $N\in \mathbb{N}$. Assume that the initial condition $z^0(\cdot)\in H^1\left (0,L;\mathbb R^3\right)$ satisfies corresponding boundary conditions. Let $\delta>0$ be a desired decay rate. Given $N$ and $\rho,\gamma>0$, assume that there exist a matrix $0\prec P\in \mathbb{R}^{3\times 3}$ and scalars $\alpha_0,\alpha_1>0$ such that $\Phi\prec 0$ (with $\Phi$ given in \eqref{eq:LMI1}) and \eqref{eq:TailSchur} hold. Then, the solution $z$ of \eqref{sys}, subject to the control law \eqref{eq:ContLaw}, \eqref{eq:GammaDef} is exponentially stable with decay rate $\delta$, meaning that \eqref{eq:stabilityineq} holds with $M$ depending polynomially on $\gamma$. In addition, the LMIs $\Phi\prec 0$ and \eqref{eq:TailSchur} are \emph{always feasible} for large enough $N$ and $\gamma$.
\end{customthm} 
\begin{IEEEproof}
	Feasibility of the LMIs \eqref{eq:LMI1} and \eqref{eq:TailSchur} implies that
	\begin{equation*}
		\dot{{V}}(t)+2\delta {V}(t)\leq 0 \Rightarrow {V}(t)\leq e^{-2\delta t}{V}(0),\ \forall t \in [0,T).
	\end{equation*}
	Thanks to the above estimate, we can extend our solution to the whole $\mathbb R_{\geq 0}$ (meaning that $T=+\infty$) and, thus, the above stability estimate holds for all $t \geq 0$. Combining the latter with \eqref{eq:Comp} yields \eqref{eq:stabilityineq}.
	
	Next, we show feasibility of the LMIs \eqref{eq:LMI1} and \eqref{eq:TailSchur} for large enough $N$ and $\gamma$. Recalling that {$Q_0+BK_0$} is Hurwitz, first choose $P\succ 0$ such that {$\operatorname{Sym}(P(Q_0+BK_0))=-I_3$}. Note that for this choice, $\bar{P} = I_N\otimes P$ in the Lyapunov functional \eqref{eq:Vfunc} satisfies $\left|\bar{P}\bar{P}^{\top} \right| =  \left|PP^{\top} \right|$ and is \emph{independent} of $N$ and $\gamma$. Setting $\rho = \gamma^{-7},$ we find that
	\begin{align}
		&\phi = -d_3 \Lambda \otimes P-\gamma \left[I_{3N}-\frac{\delta}{\gamma}I_N \otimes P \right]+\frac{\alpha_0 \ell_1}{2\gamma^7}\mathcal{M}_{N,\gamma}\nonumber\\
		&  + \frac{\alpha_1 \gamma}{2}\left(\sum_{k=1}^N\left[\left\|b_k \right\|^2 - \left| \mathcal{B}_k\right|^2 \right]  \right)\mathcal{Q}_{N\times N,\gamma}^{\top}\mathcal{Q}_{N\times N,\gamma}+\bar{P}{\overline{\mathcal{J}}_{N}}\nonumber\\
		&\overset{\eqref{Bboundedness},\eqref{Qnnbound},\eqref{eq:CalMNDef}}{\preceq}\gamma \left[-I_{3N}+\frac{\alpha_1\eta \lambda_{N+1}^{\beta}\left[\xi_{N,\gamma}+ \left| K_0\right| \right]^2}{2}I_{3N} \right.\nonumber\\
		\label{eq:phiboundUpper} &\hspace{5mm}\left.-\frac{d_3}{\gamma}\Lambda\otimes P+\frac{\alpha_0 \ell_1\sigma_N^2}{2\gamma^2}I_{3N}+\frac{\delta}{\gamma}I_N \otimes P+{\frac{1}{\gamma}\bar{P} \overline{\mathcal{J}}_{N}} \right]
	\end{align}
	where $\xi_{N,\gamma}$ and its upper bound are given in \eqref{eq:ClaGNBound} and { $\overline{\mathcal{J}}_{N}$ satisfies \eqref{eq:Jestimate}}, whence 
	\begin{equation*}
		\frac{1}{\gamma}\left| \bar{P} \overline{\mathcal{J}}_{N}\right| \leq \frac{1}{\gamma} \left| P\right| \left(\max_{1\leq i\leq 3}\left| q_{i,i} \right|+ O\left(\frac{\lambda_N}{\gamma}\right) \right)=O\left(\frac{1}{\gamma} \right). 
	\end{equation*}

	Next, consider $\Phi$ in \eqref{eq:LMI1}. Since $-\frac{\alpha_0 \rho}{2}I_{3N} = -\frac{\alpha_0}{2\gamma^7}I_{3N}\prec 0$, by Schur complement and \eqref{eq:phiboundUpper}, the LMIs \eqref{eq:LMI1} and \eqref{eq:TailSchur} are feasible if
	\begin{align}
		&-\lambda_{N+1} D +  \text{Sym}\left ( Q\right )+ \left[\frac{1}{2\alpha_0}+\frac{\alpha_0 \ell_1}{2}+\frac{1}{2\alpha_1}\right]I_3 \prec 0, \nonumber \\
		& \phi + \gamma^{-6}\bar{P}\left[\frac{2\gamma^7}{\alpha_0}I_{3N} \right]\bar{P}^{\top}\label{eq:equivcond1}\\
		&\overset{\eqref{eq:phiboundUpper}}{\prec} \gamma \left[-I_{3N}+\frac{\alpha_1\eta\lambda_{N+1}^{\beta}\left[\xi_{N,\gamma}+ \left| K_0\right| \right]^2}{2}I_{3N}+\frac{2}{\alpha_0}\bar{P}\bar{P}^{\top} \right.\nonumber \\
		&\hspace{2mm}\left.-\frac{d_3}{\gamma}\Lambda\otimes P+\frac{\alpha_0 \ell_1\sigma_N^2}{2\gamma^2}I_{3N}+\frac{\delta}{\gamma}I_N \otimes P+{\frac{1}{\gamma}\bar {P} \overline{\mathcal{J}}_{N}} \right]\prec 0.\nonumber
	\end{align}
	Recall \eqref{Bboundedness} and let $ \beta_1 \in (\beta,1)$. Set $\alpha_0 = \alpha_1^{-1} = \lambda_{N+1}^{\beta_1}$. Substituting into the first condition of \eqref{eq:equivcond1}, we have 
	\begin{equation}
		-\lambda_{N+1} D +  \text{Sym}\left ( Q\right )+ \frac{1}{2}\left[\lambda_{N+1}^{-\beta_1}+\left(\ell_1+1 \right)\lambda_{N+1}^{\beta_1}\right]I_3\prec 0. \label{eq:lambdaN+1}
	\end{equation}
	Since $\lim_{n\to + \infty}\lambda_n = +\infty$, the latter holds for all $N\geq N_*$ with $N_*$ large enough. Let $N\geq N_*$ and consider the intermediate term of the second condition in \eqref{eq:equivcond1}, given by
	\begin{align}
		& \gamma\left[\left(-1+\frac{\eta\left[\xi_{N,\gamma}+ \left| K_0\right| \right]^2}{2\lambda_{N+1}^{\beta_1-\beta}}+\frac{\lambda_{N+1}^{\beta_1} \ell\sigma_N^2}{2\gamma^2}\right)I_{3N}\right.\nonumber \\ 
		&\left.-\frac{d_3}{\gamma}\Lambda\otimes P+\frac{2}{\lambda_{N+1}^{\beta_1}}\bar{P}\bar{P}^{\top}+\frac{\delta}{\gamma}I_N \otimes P+{\frac{1}{\gamma} \bar {P}\overline{\mathcal{J}}_{N}}\right]\prec 0.\label{eq:equivcond2}
	\end{align}
	Taking into account that $\left[\xi_{N,\gamma}+ \left| K_0\right| \right]^2\leq 2\xi_{N,\gamma}^2+2\left| K_0\right|^2$, where $K_0$ is fixed, we fix $N\geq N_*$ large enough so that 
	\begin{align}
		S:=\left(-1+\frac{\eta\left| K_0\right|^2}{\lambda_{N+1}^{\beta_1-\beta}}\right)I_{3N}+\frac{2}{\lambda_{N+1}^{\beta_1}}\bar{P}\bar{P}^{\top}\prec 0. \label{eq:equivcond3}
	\end{align}
	Here, we use $0<\beta<\beta_1$ and the fact that $\left|\bar{P}\bar{P}^{\top} \right| =  \left|PP^{\top} \right|$ is \emph{independent} of $N$ and $\gamma$. Next, \eqref{eq:ClaGNBound} shows that for a fixed $N$, we have $\lim_{\gamma\to +\infty}\xi_{N,\gamma}=0$, { whereas by \eqref{eq:Jestimate}, we find that $\lim_{\gamma \to +\infty}\gamma^{-1}\left|\bar {P}\overline{\mathcal{J}}_N \right|=0$.} Hence, fixing our $N$, we have that there exists $\gamma_0>0$ such that for all $\gamma>\gamma_0$
	\begin{equation*}
		S+\left(\frac{\eta\xi_{N,\gamma}^2}{\lambda_{N+1}^{\beta_1-\beta}}+\frac{\lambda_{N+1}^{\beta_1}\ell_1 \sigma_N^2}{2\gamma^2}\right)I_{3N}+\frac{\delta}{\gamma}I_N \otimes P+{\frac{1}{\gamma}\bar {P} \overline{\mathcal{J}}_{N}}\prec 0.
	\end{equation*}
	Since $-\frac{d_3}{\gamma}\Lambda\otimes P\prec 0$ for all $\gamma>0$, we see that by choosing $\gamma>0$ large enough and taking into account \eqref{eq:equivcond3}, we obtain \eqref{eq:equivcond2}. Therefore, \eqref{eq:equivcond1} holds and the proof is concluded.
\end{IEEEproof}
\begin{customrem}{6} \it
	The proof of Theorem \ref{theorem1} allows to obtain estimates on $N$ and $\gamma$, which preserve stability of the closed-loop system, by more classical (although more conservative) arguments. Indeed, one can first find gain $K_0$, which stabilizes the matrix ${Q_0}$ and $P\succ 0$ such that $\text{Sym} \left (P({ Q_0}+BK_0) \right )\prec 0$. Second, the number of actuators $N$ can be obtained from the inequality \eqref{eq:lambdaN+1}. Given $N$, one can estimate lower bounds on $\gamma$ by guaranteeing that \eqref{eq:equivcond2} and \eqref{eq:equivcond3} hold.
\end{customrem}
\section{Numerical Example}\label{sec:example}
Consider the system \eqref{sys} with $f_2 \equiv f_3\equiv0$, $L=1$ and Dirichlet boundary conditions (i.e., $\gamma_{12}=\gamma_{22}=0$). For this case, the eigenvalues and eigenfunctions are given by {$\lambda_n = \pi^2 n^2$} and $\varphi_n(x)=\sqrt{2}\sin(\pi n x), \ n\geq 1$, respectively. We further choose $N=5$, indicator shape functions $b_j = \chi_{\left[\frac{j-1}{5},\frac{j}{5}\right)}, \quad j=1,\dots,5$ (leading to an invertible $\mathcal{B}_{N\times N}$ in \eqref{eq:BNNMat}) and 
\begin{align*}
	&D=\text{diag}\left\lbrace 2,2.5,3\right\rbrace, \quad   q_{2,1}=q_{3,2}=1, \quad {Q_1=0}.
\end{align*}
We then verify the LMIs of Theorem \ref{theorem1} with decay rate $\delta = 1 $, stabilizing gain of $Q_0$ given by
\begin{align*}
	K_0 = \begin{pmatrix}
		-3.708 & \quad -26.329 & \quad  -2.222
	\end{pmatrix}   
\end{align*}
and tuning parameters $\gamma$ and $\rho$, chosen over a grid, to obtain a Lipschitz constant $\ell_1$ which preserves feasibility of the LMIs. The LMIs were found to be feasible for $\ell_1 = 15$, with corresponding tuning parameters
\begin{align*}
	\gamma = 5,\quad \rho = 1.24\cdot 10^{-4}.
\end{align*}
Note that the value of the Lipschitz constant for which feasibility holds is large. This numerical simulation demonstrates the efficiency of the proposed approach.

\section{Conclusion}\label{sec:conclusion}



A novel method for exponential stabilization of a system of parabolic PDEs interconnected via a reaction matrix and a Lipschitz nonlinearity was presented. The proposed approach is based on three steps: modal decomposition, transformation of the mode ODEs into an equivalent target system, and design of a high-gain state-feedback controller that compensates the nonlinearity. Future research topics may include other classes of parabolic PDEs and more general interconnections. 

\bibliographystyle{IEEEtran}
\bibliography{IEEEabrv,Bibliography111022}
\end{document}